\documentclass[12pt]{article}
\usepackage{amsmath}
\usepackage{amssymb}
\usepackage{amsfonts}
\usepackage{amsthm}
\usepackage{enumerate}

\def\init{\setcounter{equation}{0}}

\newtheorem{theorem}{Theorem}[section]
\newtheorem{proposition}[theorem]{Proposition}

\newtheorem{definition}[theorem]{Definition}
\newtheorem{corollary}[theorem]{Corollary}

\def\Li{{\operatorname{Li}}}
\def\Res{{\operatorname{Res}}}
\def\m{{\operatorname{m}}}

\def\d{{\operatorname{d}}}

\def\Re{\hbox{Re}\,}
\def\Im{\hbox{Im}\,}

\def\mod{\hbox{mod}\,}
\def\p{\partial}

\def\ii{{\mathrm{i}}}
\def\e{{\mathrm{e}}}
\def\r{{\, \texttt{r}}}

\numberwithin{equation}{section}
\newenvironment{acknowledgements}{\noindent{\bf Acknowledgements}\bigskip}{}

\begin{document}

\title{Functional equations for Mahler measures of genus-one curves}
\author{Matilde N. Lal\'{\i}n, Mathew D. Rogers \\
        \small{\textit{Department of Mathematics, University of British
        Columbia}}\\
        \small{\textit{Vancouver, BC, V6T-1Z2, Canada}}}

\maketitle

\abstract{In this paper we will establish functional equations for
 Mahler measures of families of  genus-one
two-variable polynomials. These families were previously studied by Beauville, and their Mahler measures were considered
by Boyd, Rodriguez-Villegas,
Bertin, Zagier, and Stienstra.
Our functional equations allow us to prove identities between Mahler measures that were conjectured by Boyd.  As a
corollary, we also establish some new transformations for
hypergeometric functions.}

\medskip


\section{History and introduction}
\label{intro} \init

    The goal of this paper is to establish identities between the logarithmic
Mahler measures of polynomials with zero varieties corresponding to
genus-one curves.  Recall that the logarithmic Mahler measure (which
we shall henceforth simply refer to as the Mahler measure) of an
$n$-variable Laurent polynomial $P(x_1,x_2,\dots,x_n)$ is defined by
\begin{equation*}
\m\left(P(x_1,\dots,x_n)\right)=\int_{0}^{1}\dots\int_{0}^{1}\log\big\vert
P\left(\e^{2\pi \ii\theta_1},\dots,\e^{2\pi \ii\theta_n}\right)\big\vert
\d\theta_1\dots\d\theta_n.
\end{equation*}
Many difficult questions surround the special functions defined by
Mahler measures of elliptic curves.

The first example of the Mahler measure of a genus-one curve was
studied by Boyd \cite{Bo1} and Deninger \cite{D}. Boyd found that
\begin{equation}\label{eq:bd}
\m\left(1+x+\frac{1}{x}+y+\frac{1}{y}\right)\stackrel{?}{=} L'(E,0),
\end{equation}
where $E$ denotes the elliptic curve of conductor 15 that is the
projective closure of $1+x+\frac{1}{x}+y+\frac{1}{y}=0$.  As usual,
$L(E,s)$ is its $L$-function, and the question mark above the equals
sign indicates numerical equality verified up to 28 decimal places.

Deninger \cite{D} gave an interesting interpretation of this
formula. He obtained the Mahler measure  by evaluating the
Bloch regulator of an element $\{x,y\}$ from a certain $K$-group. In
other words, the Mahler measure is given by a value of an
Eisenstein-Kronecker series. Therefore Bloch's and Beilinson's
conjectures predict that
\[\m\left(1+x+\frac{1}{x}+y+\frac{1}{y}\right)= cL'(E,0),\]
where $c$ is some rational number. Let us add that, even if
Beilinson's conjectures were known to be true, this would not
suffice to prove equality \eqref{eq:bd}, since we still would not
know the height of the rational number $c$.

This picture applies to other situations as well.  Boyd \cite{Bo1}
performed extensive numerical computations within the family of
polynomials $k+x+\frac{1}{x}+y+\frac{1}{y}$, as well as within some
other genus-one families. Boyd's numerical searches led him to
conjecture identities such as
\[\m\left(5+x+\frac{1}{x}+y+\frac{1}{y}\right)\stackrel{?}{=} 6 \m\left(1+x+\frac{1}{x}+y+\frac{1}{y}\right),\]
\[\m\left(8+x+\frac{1}{x}+y+\frac{1}{y}\right)\stackrel{?}{=} 4 \m\left(2+x+\frac{1}{x}+y+\frac{1}{y}\right).\]
Boyd conjectured conditions predicting when formulas like Eq.
\eqref{eq:bd} should exist for the Mahler measures of polynomials
with integral coefficients. This was further studied by
Rodriguez-Villegas \cite{RV} who interpreted these
conditions in the context of Bloch's and Beilinson's conjectures. Furthermore, Rodriguez-Villegas used modular forms to
express the Mahler measures as Kronecker-Eisenstein series in more
general cases. In turn, this allowed him to prove some equalities
such as
\begin{equation} \label{eq:frv1}
\m\left(4\sqrt{2} +x +\frac{1}{x}+y +\frac{1}{y}\right)=L'\left(E_{4\sqrt{2}},0\right),
\end{equation}
\begin{equation}\label{eq:frv2}
\m\left(3 \sqrt{2} +x +\frac{1}{x}+y +\frac{1}{y}\right)=\frac{5}{2} L'\left(E_{3\sqrt{2}},0\right).
\end{equation}
The first equality can be proved using the fact that the
corresponding elliptic curve has complex multiplication, and
therefore the conjectures are known for this case due to Bloch
\cite{B}. The second equality depends on the fact that one has the
modular curve $X_0(24)$, and the conjectures then follow from a
result of Beilinson.

Rodriguez-Villegas \cite{RV2} subsequently used the relationship
between Mahler measures and regulators to prove a conjecture of Boyd
\cite{Bo1}:
\[\m(y^2+2xy+y-x^3-2x^2-x)=\frac{5}{7}\m(y^2+4xy+y-x^3+x^2).\]
It is important to point out that he proved this identity without
actually expressing the Mahler measures in terms of $L$-series.
Bertin \cite{bert1} has also proved similar identities using these
ideas.

    Although the conjecture in Eq. \eqref{eq:bd} remains open, we will in fact prove two of Boyd's other conjectures this paper.
\begin{theorem}\label{big theorem} The following identities are true:
\begin{align}
m\left(2+x+\frac{1}{x}+y+\frac{1}{y}\right)=&L'\left(E_{3\sqrt{2}},0\right),\\
m\left(8+x+\frac{1}{x}+y+\frac{1}{y}\right)=&4L'\left(E_{3\sqrt{2}},0\right).
\end{align}
\end{theorem}\noindent
Our proof of Theorem \ref{big theorem} follows from combining two
interesting ``functional equations" for the function
\[m(k):=\m\left(k+x+\frac{1}{x}+y+\frac{1}{y}\right).\]
Kurokawa and Ochiai \cite{KO} recently proved the first functional
equation.  They showed that if $k\in\mathbb{R}\backslash \{0\}$:
\begin{equation} \label{eq:ko}
m\left(4k^2\right)+m\left(\frac{4}{k^2}\right) = 2
m\left(2\left(k+\frac{1}{k}\right)\right).
\end{equation}
In Section \ref{regul} we use regulators to give a new proof of Eq.
\eqref{eq:ko}.  We will also prove a second functional equation in
Section \ref{modular equations subsection} using $q$-series.  In
particular, if $k$ is nonzero and $|k|<1$:
\begin{equation} \label{eq:first}
m\left(2\left(k+\frac{1}{k}\right)\right)+m\left(2\left(\ii
k+\frac{1}{\ii k}\right)\right)=  m\left(\frac{4}{k^2}\right).
\end{equation}
Theorem \ref{big theorem} follows from setting $k=1/\sqrt{2}$ in
both identities, and then showing that
$5m\left(\ii\sqrt{2}\right)=3m\left(3\sqrt{2}\right)$.  We have
proved this final equality in Section \ref{sec:sqrt{2}}.

This paper is divided into two sections of roughly equal length. In
Section \ref{q series summary} we will prove more identities like
Eq. \eqref{eq:first}, which arise from expanding Mahler measures in
$q$-series.  In particular, we will look at identities for four
special functions defined by the Mahler measures of genus-one curves
(see equations \eqref{eq:m} through \eqref{eq:h} for notation).
Equation \eqref{Hecke eigenvalue equation} undoubtedly constitutes
the most important result in this part of the paper, since it
implies that infinitely many identities like Eq. \eqref{eq:first}
exist. Subsections \ref{modular equations subsection} and
\ref{higher modular equations subsection} are mostly devoted to
transforming special cases of Eq. \eqref{Hecke eigenvalue equation}
into interesting identities between the Mahler measures of rational
polynomials.  While the theorems in those subsections rely heavily
on Ramanujan's theory of modular equations to alternative bases, we
have attempted to maximize readability by eliminating $q$-series
manipulation wherever possible.  Finally, we have devoted Subsection
\ref{computational subsection} to proving some useful computational
formulas.  As a corollary we establish several new transformations
for  hypergeometric functions, including:
\begin{equation}\label{picard fuchs}
\begin{split}
\sum_{n=0}^{\infty}&\left(\frac{k(1-k)^2}{(1+k)^2}\right)^{n}\sum_{j=0}^{n}{n\choose
j}^2{n+j\choose
j}\\
&=\frac{(1+k)^2}{\sqrt{\left(1+k^2\right)\left(\left(1-k-k^2\right)^2-5k^2\right)}}{_2F_1}\left(\frac{1}{4},\frac{3}{4};1;
\frac{64k^5\left(1+k-k^2\right)}{\left(1+k^2\right)^2\left(\left(1-k-k^2\right)^2-5k^2\right)^2}\right).
\end{split}
\end{equation}

We have devoted Section \ref{regul} to further studying the
relationship between Mahler measures and regulators. We show how to
recover the Mahler measure $q$-series expansions and the
Kronecker-Einsenstein series directly from Bloch's formula for the
regulator.  This in turn shows that the Mahler measure identities
can be viewed as consequences of functional identities for the
elliptic dilogarithm.

  Many of the identities in this paper can be interpreted from
both a regulator perspective, and from a $q$-series perspective. The
advantage of the $q$-series approach is that it simplifies the
process of finding new identities.  The fundamental result in
Section \ref{q series summary}, Eq. \eqref{Hecke eigenvalue
equation}, follows easily from the Mahler measure $q$-series
expansions.  Unfortunately the $q$-series approach does not provide
an easy way to explain identities like Eq. \eqref{eq:ko}. Unlike
most of the other formulas in Section \ref{q series summary},
Kurokawa's and Ochiai's result \textit{does not} follow from Eq.
\eqref{Hecke eigenvalue equation}. An advantage of the regulator
approach, is that it enables us to construct proofs of both Eq.
\eqref{eq:ko} and Eq. \eqref{eq:first} from a unified perspective.
Additionally, the regulator approach seem to provide the only way to
prove the final step in Theorem \ref{big theorem}, namely to show
that $5m\left(\ii\sqrt{2}\right)=3m\left(3\sqrt{2}\right)$.  Thus, a
complete view of this subject matter should incorporate both
regulator and $q$-series perspectives.

\section{Mahler measures and $q$-series}
\label{q series summary} \init
In this paper we will consider four important functions defined by
Mahler measures:
\begin{align}
\mu(t)=&\m\left(\frac{4}{\sqrt{t}}+x+\frac{1}{x}+y+\frac{1}{y}\right), \label{eq:m}\\
n(t)=&\m\left(x^3+y^3+1-\frac{3}{t^{1/3}}x y\right), \label{eq:n} \\
g(t)=&\m\left((x+y)(x+1)(y+1)-\frac{1}{t}x y\right), \label{eq:g}\\
r(t)=&\m\left((x+y+1)(x+1)(y+1)-\frac{1}{t} x y\right). \label{eq:h}
\end{align}
Throughout Section \ref{q series summary} will use the notation
$\mu(t)=m\left(\frac{4}{\sqrt{t}}\right)$ for convenience.  Recall
from \cite{RV} and \cite{St}, that each of these functions has a
simple $q$-series expansion when $t$ is parameterized correctly.  To
summarize, if we let $(x;q)_{\infty}=(1-x)(1-x q)\left(1-x
q^2\right)\dots$, and
\begin{align}
M(q)&=16 q\frac{(q;q)^8_{\infty}\left(q^4;q^4\right)^{16}_\infty}{\left(q^2;q^2\right)^{24}_{\infty}},\label{M(q)}\\
N(q)&=\frac{27q\left(q^3;q^3\right)^{12}_{\infty}}{(q;q)^{12}_{\infty}+27q\left(q^3;q^3\right)^{12}_{\infty}},\label{N(q)}\\
G(q)&=q^{1/3}\frac{\left(q;q^2\right)_{\infty}}{\left(q^3;q^6\right)^3_{\infty}},\label{G(q)}\\
R(q)&=q^{1/5}\frac{\left(q;q^5\right)_{\infty}\left(q^4;q^5\right)_{\infty}}{\left(q^2;q^5\right)_{\infty}\left(q^3;q^5\right)_{\infty}},\label{R(q)}
\end{align}
then for $|q|$ sufficiently small
\begin{align}
\mu\left(M(q)\right)=&-\Re\left[\frac{1}{2}\log(q)
+2\sum_{j=1}^{\infty}j\chi_{-4}(j)\log\left(1-q^j\right)\right],\label{m(t) q series}\\
n\left(N(q)\right)=&-\Re\left[\frac{1}{3}\log(q)+3\sum_{j=1}^{\infty}j\chi_{-3}(j)\log(1-q^j)\right],\label{n(t) q series}\\
g\left(G^3(q)\right)=&-\Re\left[\log(q)+\sum_{j=1}^{\infty}(-1)^{j-1}j\chi_{-3}(j)\log(1-q^j)\right],\label{g(t) q series}\\
r\left(R^5(q)\right)=&-\Re\left[\log(q)+
\sum_{j=1}^{\infty}j\Re\left[(2-\ii)\chi_r(j)\right]\log\left(1-q^j\right)\right].\label{r(t) q series}
\end{align}
In particular, $\chi_{-3}(j)$ and $\chi_{-4}(j)$ are the usual
Dirichlet characters, and $\chi_{r}(j)$ is the character of
conductor five with $\chi_{r}(2)=\ii$.  We have used the notation
$G(q)$ and $R(q)$, as opposed to something like
$\tilde{G}(q)=G^3(q)$, in order to preserve Ramanujan's notation.
As usual, $G(q)$ corresponds to Ramanujan's cubic continued
fraction, and $R(q)$ corresponds to the Rogers-Ramanujan continued
fraction \cite{BA}.

    The first important application of the $q$-series expansions
is that they can be used to calculate the Mahler measures
numerically. For example, we can calculate $\mu\left(1/10\right)$
with Eq. \eqref{m(t) q series}, provided that we can first determine
a value of $q$ for which $M(q)=1/10$. Fortunately, the theory of
elliptic functions shows that if $ \alpha=M(q)$, then
\begin{equation}\label{elliptic nome}
q=\exp\left(-\pi\frac{{_2F_1}\left(\frac{1}{2},\frac{1}{2};1;1-\alpha\right)}
{{_2F_1}\left(\frac{1}{2},\frac{1}{2};1;\alpha\right)}\right).
\end{equation}
Using Eq. \eqref{elliptic nome} we easily compute $q=.01975\dots$,
and it follows that $\mu(1/10)=2.524718\dots$  The function defined
in Eq. \eqref{elliptic nome} is called the \textit{elliptic nome},
and is sometimes denoted by $q_2(\alpha)$.  Theorem \ref{inversion
theorem} provides similarly explicit inversion formulas for Eqs.
\eqref{M(q)} through \eqref{R(q)}.

The second, and perhaps more significant fact that follows from
these $q$-series, is that linear dependencies exist between the
Mahler measures. In particular, if
\[f(q)\in\left\{\mu\left(M(q)\right),n\left(N(q)\right),g\left(G^3(q)\right),r\left(R^5(q)\right)\right\},\]
then for an appropriate prime $p$
\begin{equation}\label{Hecke eigenvalue equation}
\sum_{j=0}^{p-1}f\left(\e^{2\pi\ii j/p}
q\right)=(1+p^2\chi(p))f\left(q^p\right)-p\chi(p)
f\left(q^{p^2}\right),
\end{equation}
where $\chi(j)$ is the character from the relevant $q$-series. The
prime $p$ satisfies the restriction that $p\not=2$ when
$f(q)=g\left(G^3(q)\right)$, and $p\not\equiv 2,3\left(\mod
5\right)$ when $f(q)=r\left(R^5(q)\right)$.  The astute reader will
immediately recognize that Eq. \eqref{Hecke eigenvalue equation} is
essentially a Hecke eigenvalue equation.  A careful analysis of the
exceptional case that occurs when $p=2$ and
$f(q)=g\left(G^3(q)\right)$ leads to the important and surprising
inverse relation:
\begin{equation}\label{functional equation g and h mixed}
\begin{split}
3 n\left(N(q)\right)&=g\left(G^3\left(q\right)\right)-8g\left(G^3\left(-q\right)\right)+4g\left(G^3\left(q^2\right)\right),\\
3
g\left(G^3(q)\right)&=n\left(N(q)\right)+4n\left(N\left(q^2\right)\right).
\end{split}
\end{equation}
In the next two subsections we will discuss methods for transforming
Eq. \eqref{Hecke eigenvalue equation} and Eq. \eqref{functional
equation g and h mixed} into so-called functional equations.

\subsection{Functional equations from modular equations}
\label{modular equations subsection}

    Since the primary goal of this paper is to find relations
between the Mahler measures of \textit{rational} (or at least
algebraic) polynomials, we will require modular equations to
simplify our results.  For example, consider Eq. \eqref{Hecke
eigenvalue equation} when $f(q)=\mu\left(M(q)\right)$ and $p=2$:
\begin{equation}\label{m(t) second order pre equation}
\mu\left(M(q)\right)+\mu\left(M(-q)\right)=\mu\left(M\left(q^2\right)\right).
\end{equation}
For our purposes, Eq. \eqref{m(t) second order pre equation} is only
interesting if $M(q)$, $M(-q)$, and $M\left(q^2\right)$ are all
simultaneously algebraic.  Fortunately, it turns out that $M(q)$ and
$M\left(q^2\right)$ (hence also $M(-q)$ and $M\left(q^2\right)$)
satisfy a well known polynomial relation.

\begin{definition}Suppose that
$F(q)\in\left\{M(q),N(q),G(q),R(q)\right\}$.  An $n$'th degree
modular equation is an algebraic relation between $F(q)$ and
$F(q^n)$.
\end{definition}

We will not need to derive any new modular equations in this paper.
Berndt proved virtually all of the necessary modular equations while
editing Ramanujan's notebooks, see \cite{BA}, \cite{Be2},
\cite{Be3}, and \cite{Be5}. Ramanujan seems to have arrived at most
of his modular equations through complicated $q$-series
manipulations (of course this is speculation since he did not write
down any proofs!). Modular equations involving $M(q)$ correspond to
the classical modular equations \cite{Be3},  relations for $N(q)$
correspond to Ramanujan's signature three modular equations
\cite{Be5}, and most of the known modular equations for $G(q)$ and
$R(q)$ appear in \cite{BA}.

Now we can finish simplifying Eq. \eqref{m(t) second order pre
equation}.  Since the classical second-degree modular equation shows
that whenever $|q|<1$
\begin{equation*}
\frac{4M\left(q^2\right)}{\left(1+M\left(q^2\right)\right)^2}=\left(\frac{M(q)}{M(q)-2}\right)^2,
\end{equation*}
we easily obtain the parameterizations:
$M(q)=\frac{4k^2}{\left(1+k^2\right)^2}$,
$M(-q)=\frac{-4k^2}{\left(1-k^2\right)^2}$, and
$M\left(q^2\right)=k^4$.  Substituting these parametric formulas
into Eq. \eqref{m(t) second order pre equation} yields:
\begin{theorem}\label{theorem on functional equation p2 M(q)} The following identity holds whenever $|k|<1$:
\begin{equation}\label{functional equation p2 M(q)}
\begin{split}
\m\left(\frac{4}{k^2}+x+\frac{1}{x}+y+\frac{1}{y}\right)=&\m\left(2\left(k+\frac{1}{k}\right)+x+\frac{1}{x}+y+\frac{1}{y}\right)\\
&+\m\left(2\ii\left(k-\frac{1}{k}\right)+x+\frac{1}{x}+y+\frac{1}{y}\right).
\end{split}
\end{equation}
\end{theorem}



    We need to make a few remarks about working with modular equations
before proving the main theorem in this section.  Suppose that for
some algebraic function $P(X,Y)$:
\begin{equation*}
P\left(F(q),F\left(q^p\right)\right)=0,
\end{equation*}
where $F(q)\in\left\{M(q),N(q),G(q),R(q)\right\}$.  Using the
elementary change of variables, $q\rightarrow \e^{2\pi\ii j/p} q$,
it follows that $P\left(F\left(\e^{2\pi \ii j/p}
q\right),F\left(q^p\right)\right)=0$ for every
$j\in\{0,1,\dots,p-1\}$.  If $P(X,Y)$ is symmetric in $X$ and $Y$,
it also follows that
$P\left(F\left(q^{p^2}\right),F\left(q^p\right)\right)=0$.
Therefore, if $P(X,Y)$ is sufficiently simple (for example a
symmetric genus-zero polynomial), we can find simultaneous
parameterizations for $F\left(q^p\right)$, $F\left(q^{p^2}\right)$,
and $F\left(\e^{2\pi\ii j/p} q\right)$ for all $j$. In such an
instance, Eq. \eqref{Hecke eigenvalue equation} reduces to an
interesting functional equation for one of the four Mahler measures
$\{\mu(t),n(t),g(t),r(t)\}$.  Five basic functional equations follow
from applying these ideas to Eq. \eqref{Hecke eigenvalue equation}.

\begin{theorem}\label{functional equations theorem}
For $|k|<1$ and $k\not=0$, we have
\begin{equation}\label{functional equation p=2 m(t)}
\mu\left(\frac{4k^2}{\left(1+k^2\right)^2}\right)+\mu\left(\frac{-4k^2}{\left(1-k^2\right)^2}\right)=\mu\left(k^4\right).
\end{equation}
The following identities hold for $|u|$ sufficiently small but
non-zero:
\begin{equation}\label{functional equation p=2 n(t)}
\begin{split}
n\left(\frac{27u(1+u)^4}{2\left(1+4u+u^2\right)^3}\right)+&n\left(-\frac{27u(1+u)}{2\left(1-2u-2u^2\right)^3}\right)\\
&=2n\left(\frac{27u^4(1+u)}{2\left(2+2u-u^2\right)^3}\right)
-3n\left(\frac{27u^2(1+u)^2}{4\left(1+u+u^2\right)^3}\right).
\end{split}
\end{equation}
If $\zeta_3=\e^{2\pi\ii/3}$, and
$Y(t)=1-\left(\frac{1-t}{1+2t}\right)^3$, then
\begin{equation}\label{functional equation p=3 n(t)}
n\left(u^3\right)=\sum_{j=0}^{2}n\left(Y\left(\zeta_3^j
u\right)\right).
\end{equation}
If $\zeta_3=\e^{2\pi\ii/3}$, and
$Y(t)=t\left(\frac{1-t+t^2}{1+2t+4t^2}\right)$, then
\begin{equation}\label{functional equation p=3 g(t)}
g(u^3)=\sum_{j=0}^{2}g\left(Y\left(\zeta_3^j u\right)\right).
\end{equation}
If $\zeta_5=\e^{2\pi \ii /5}$, and
$Y(t)=t\left(\frac{1-2t+4t^2-3t^3+t^4}{1+3t+4t^2+2t^3+t^4}\right)$,
then
\begin{equation}\label{functional equation p=5 r(t)}
r\left(u^5\right)=\sum_{j=0}^{4}r\left(Y\left(\zeta_5^j
u\right)\right).
\end{equation}
\end{theorem}
\begin{proof}
We have already sketched a proof of Eq. \eqref{functional equation
p=2 m(t)} in the discussion preceding Theorem \ref{theorem on
functional equation p2 M(q)}.

    The proof of Eq. \eqref{functional equation p=2 n(t)} requires
the second-degree modular equation from Ramanujan's theory of
signature three.  If $\beta=N\left(q^2\right)$, and
$\alpha\in\left\{N(q),N(-q),N\left(q^4\right)\right\}$, then
\begin{equation}\label{modular equation sig 3 deg 2}
27\alpha\beta(1-\alpha)(1-\beta)-\left(\alpha+\beta-2\alpha\beta\right)^3=0.
\end{equation}
If we choose $u$ so that
$N\left(q^2\right)=\frac{27u^2(1+u)^2}{4\left(1+u+u^2\right)^3}$,
then we can use Eq. \eqref{modular equation sig 3 deg 2} to easily
verify that $N(q)=\frac{27u(1+u)^4}{2\left(1+4u+u^2\right)^3}$,
$N(-q)=-\frac{27u(1+u)}{2\left(1-2u-2u^2\right)^3}$, and
$N\left(q^4\right)=\frac{27u^4(1+u)}{2\left(2+2u-u^2\right)^3}$. The
proof of Eq. \eqref{functional equation p=2 n(t)} follows from
applying these parameterizations to Eq. \eqref{Hecke eigenvalue
equation} when $f(q)=n\left(N(q)\right)$, and $p=2$.

    The proof of Eq. \eqref{functional equation p=3 n(t)} requires
Ramanujan's third-degree, signature three modular equation.  In
particular, if $\alpha=N(q)$ and $\beta=N(q^3)$, then
\begin{equation}\label{modular equation sig 3 deg 3}
\alpha=1-\left(\frac{1-\beta^{1/3}}{1+2\beta^{1/3}}\right)^3
=Y\left(\beta^{1/3}\right).
\end{equation}
Since $N^{1/3}(q^3)=q\times\{\text{power series in $q^3$}\}$, a
short computation shows that $N(\zeta_3^j
q)=Y\left(\zeta_3^{j}N^{1/3}\left(q^3\right)\right)$ for all
$j\in\{0,1,2\}$. Choosing $u$ such that $N(q^3)=u^3$, we must have
$N\left(\zeta_3^{j} q\right)=Y\left(\zeta_3^{j} u\right)$.  Eq.
\eqref{functional equation p=3 n(t)} follows from applying these
parametric formulas to Eq. \eqref{Hecke eigenvalue equation} when
$f(q)=n\left(N(q)\right)$, and $p=3$.

    Since the proofs of equations \eqref{functional equation p=3 g(t)} and \eqref{functional equation p=5
r(t)} rely on similar arguments to the proof of Eq.
\eqref{functional equation p=3 n(t)}, we will simply state the
prerequisite modular equations. In particular, Eq. \eqref{functional
equation p=3 g(t)} follows from Ramanujan's third-degree modular
equation for the cubic continued fraction.  If $\alpha=G(q)$ and
$\beta=G\left(q^3\right)$, then
\begin{equation}
\alpha^3=\beta\left(\frac{1-\beta+\beta^2}{1+2\beta+4\beta^2}\right).
\end{equation}
Similarly, Eq. \eqref{functional equation p=5 r(t)} follows from the
fifth-degree modular equation for the Rogers-Ramanujan continued
fraction.  In particular, if $\alpha=R(q)$ and
$\beta=R\left(q^5\right)$
\begin{equation}
\alpha^5=\beta\left(\frac{1-2\beta+4\beta^2-3\beta^3+\beta^4}{1+3\beta+4\beta^2+2\beta^3+\beta^4}\right).
\end{equation}$\blacksquare$
\end{proof}
    The functional equations in Theorem \ref{functional equations
    theorem} only hold in restricted subsets of $\mathbb{C}$.  To explain this phenomenon
we will go back to Eq. \eqref{Hecke eigenvalue equation}. As a
general rule, we have to restrict $q$ to values for which
\textit{none} of the Mahler measure integrals in Eq. \eqref{Hecke
eigenvalue equation} vanish on the unit torus. In other words, we
can only consider the set of $q$'s for which each term in Eq.
\eqref{Hecke eigenvalue equation} can be calculated from the
appropriate $q$-series.  Next, we may need to further restrict the
domain of $q$ depending on where the relevant parametric formulas
hold. For example, parameterizations such as
$N(q)=\frac{27u(1+u)^4}{2(1+4u+u^2)^3}$ and
$N(q^2)=\frac{27u^2(1+u)^2}{4(1+u+u^2)^3}$ hold for $|q|$
sufficiently small, but fail when $q$ is close to $1$. After
determining the domain of $q$, we can calculate the domain of $u$ by
solving a parametric equation to express $u$ in terms of a
$q$-series.

\begin{theorem}\label{functional equation n and g theorem} For $|p|$ sufficiently
small but non-zero
\begin{equation}\label{functional equation g in terms of n}
3g(p)=n\left(\frac{27p}{(1+4p)^3}\right)+4n\left(\frac{27
p^2}{(1-2p)^3}\right).
\end{equation}
Furthermore, for $|u|$ sufficiently small but non-zero
\begin{equation}\label{functional equation n in terms of g}
\begin{split}
3n\left(\frac{27u(1+u)^4}{2(1+4u+u^2)^3}\right)=&g\left(\frac{u}{2(1+u)^2}\right)
-8g\left(-\frac{u(1+u)}{2}\right)\\
&+4g\left(\frac{u^2}{4(1+u)}\right).
\end{split}
\end{equation}
\end{theorem}
\begin{proof}We will prove Eq. \eqref{functional equation n in terms of g} first.
Recall that Eq. \eqref{functional equation g and h mixed} shows that
\begin{equation*}
3n\left(N(q)\right)=g\left(G^3(q)\right)-8g\left(G^3(-q)\right)+4g\left(G^3\left(q^2\right)\right).
\end{equation*}
Let us suppose that $q=q_2\left(\frac{u(2+u)^3}{(1+2u)^3}\right)$,
where $q_2(\alpha)$ is the elliptic nome.  Classical eta function
inversion formulas (which we shall omit here) show that for $|u|$
sufficiently small: $G^3(q)=\frac{u}{2(1+u)^2}$,
$G^3(-q)=-\frac{u(1+u)}{2}$,
$G^3\left(q^2\right)=\frac{u^2}{4(1+u)}$,
$N(q)=\frac{27u(1+u)^4}{2\left(1+4u+u^2\right)^3}$, and
$N\left(q^2\right)=\frac{27u^2(1+u)^2}{4\left(1+u+u^2\right)^3}$.

    To prove Eq. \eqref{functional equation g in terms of n} first recall
recall that
\begin{equation*}
3
g\left(G^3(q)\right)=n\left(N(q)\right)+4n\left(N\left(q^2\right)\right).
\end{equation*}
If we let $p=\frac{u}{2(1+u)^2}$, then it follows that $G^3(q)=p$,
$N(q)=\frac{27p}{(1+4p)^3}$, and
$N\left(q^2\right)=\frac{27p^2}{(1+2p)^3}$.$\blacksquare$
\end{proof}
Theorem \ref{functional equation n and g theorem} shows that $g(t)$
and $n(t)$ are essentially interchangeable.  In Section
\ref{computational subsection} we will use Eq. \eqref{functional
equation g in terms of n} to derive an extremely useful formula for
calculating $g(t)$ numerically.
\subsection{Identities arising from higher modular equations}
\label{higher modular equations subsection}
The seven functional equations presented in Section \ref{modular
equations subsection} are certainly not the only interesting
formulas that follow from Eq. \eqref{Hecke eigenvalue equation}.
Rather those results represent the subset of functional equations in
which every Mahler measure depends on a rational argument (possibly
in a cyclotomic field).  If we consider the higher modular
equations, then we can establish formulas involving the Mahler
measures of the modular polynomials themselves.  Eq. \eqref{mahler
measure of the modular polynomial p=3} is the simplest formula in
this class of results. 

    Consider Eq. \eqref{Hecke eigenvalue equation} when
$p=3$ and $f(q)=\mu\left(M(q)\right)$:
\begin{equation}\label{functional equation m(M(q)) p=3}
\sum_{j=0}^{2}\mu\left(M\left(\zeta_3^j q\right)\right)=-8
\mu\left(M\left(q^3\right)\right)+3
\mu\left(M\left(q^{9}\right)\right).
\end{equation}
The third-degree modular equation shows that if
$\alpha\in\left\{M\left(q\right),M\left(\zeta_3
q\right),M\left(\zeta_3^2 q\right),M\left(q^9\right)\right\}$, and
$\beta=M\left(q^3\right)$, then
\begin{equation}\label{modular equation p=3 f(q)=M(q)}
G_3(\alpha,\beta):=(\alpha^2+\beta^2+6\alpha\beta)^2-16\alpha\beta\left(4(1+\alpha\beta)-3(\alpha+\beta)\right)^2=0.
\end{equation}
Since $G_3(\alpha,\beta)=0$ defines a curve with genus greater than
zero, it is impossible to find simultaneous rational
parameterizations for all four zeros in $\alpha$.  For example, if
we let $\beta=M(q^3)=p(2+p)^3/(1+2p)^3$, then we can obtain the
rational expression $M\left(q^9\right)=p^3(2+p)/(1+2p)$, and three
messy formulas involving radicals for the other zeros. Despite this
difficulty, Eq. \eqref{functional equation m(M(q)) p=3} still
reduces to an interesting formula if we recall the factorization
\begin{equation}\label{G3 factored}
G_3\left(\alpha,M\left(q^3\right)\right)=\left(\alpha-M\left(q^9\right)\right)\prod_{j=0}^{2}\left(\alpha-M\left(\zeta_3^j
q\right)\right),
\end{equation}
and then use the fact that Mahler measure satisfies
$\m(P)+\m(Q)=\m\left(PQ\right)$.

\begin{theorem} If $G_3(\alpha,\beta)$ is defined in Eq.
\eqref{modular equation p=3 f(q)=M(q)}, then for $|p|$ sufficiently
small but non-zero
\begin{equation}\label{mahler measure of the modular polynomial p=3}
\begin{split}
\m&\left(G_3\left(\frac{\left(x+x^{-1}\right)^2\left(y+y^{-1}\right)^2}{16},\frac{1}{p}\left(\frac{1+2p}{2+p}\right)^3\right)\right)\\
&=-16\log(2)-16\mu\left(p\left(\frac{2+p}{1+2p}\right)^3\right)+8\mu\left(p^3\left(\frac{2+p}{1+2p}\right)\right).
\end{split}
\end{equation}
\end{theorem}
\begin{proof}
First notice that from the elementary properties of Mahler's measure
\begin{equation*}
\mu(t)=\frac{1}{2}\m\left(\frac{16}{\left(x+x^{-1}\right)^2\left(y+y^{-1}\right)^2}-t\right)-\frac{1}{2}\log|t|.
\end{equation*}
Applying this identity to Eq. \eqref{functional equation m(M(q))
p=3}, and then appealing to Eq. \eqref{G3 factored} yields
\begin{equation*}
\begin{split}
\m\left(G_3\left(\frac{16}{\left(x+x^{-1}\right)^2\left(y+y^{-1}\right)^2},M\left(q^3\right)\right)\right)=&\log\left|M\left(q\right)M\left(\zeta_3
q\right)M\left(\zeta_3^2 q\right)M\left(q^9\right)\right|\\
&-16
\mu\left(M\left(q^3\right)\right)+8\mu\left(M\left(q^{9}\right)\right).
\end{split}
\end{equation*}
Elementary $q$-product manipulations show that $
M^4\left(q^3\right)=M\left(q\right)M\left(\zeta_3
q\right)M\left(\zeta_3^2 q\right)M\left(q^9\right)$, and since
$\alpha^4\beta^4
G_3\left(\frac{1}{\alpha},\frac{1}{\beta}\right)=G_3(\alpha,\beta)$,
we obtain
\begin{equation*}
\begin{split}
\m\left(G_3\left(\frac{\left(x+x^{-1}\right)^2\left(y+y^{-1}\right)^2}{16},\frac{1}{M\left(q^3\right)}\right)\right)=&-16\log(2)
-16 \mu\left(M\left(q^3\right)\right)\\&+8
\mu\left(M\left(q^{9}\right)\right).
\end{split}
\end{equation*}
Finally, if we choose $p$ so that
$M\left(q^3\right)=p\left(\frac{2+p}{1+2p}\right)^3$, then
$M\left(q^9\right)=p^3\left(\frac{2+p}{1+2p}\right)$, and the
theorem follows. $\blacksquare$
\end{proof}

Although we completely eliminated the $q$-series expressions from
Eq. \eqref{mahler measure of the modular polynomial p=3}, this is
not necessarily desirable (or even possible) in more complicated
examples.  Consider the identity involving resultants which follows
from Eq. \eqref{Hecke eigenvalue equation} (and some manipulation)
when $p=11$ and $f(q)=r\left(R^5(q)\right)$:
\begin{equation}\label{mahler measure for p=11 r(t)}
\begin{split}
\m&\left(\mathop{\Res}_z\left[z^5-\frac{x
y}{(x+1)(y+1)(x+y+1)},P\left(z,R^5\left(q\right)\right)\right]\right)\\
&=-12\m\left(1+x+y\right)+12\log\left|R^5\left(q\right)\right|+122r\left(R^5\left(q\right)\right)-11r\left(R^5\left(q^{11}\right)\right).
\end{split}
\end{equation}
In this formula $P(u,v)$ is the polynomial
\[P(u,v)=uv(1-11v^5-v^{10})(1-11u^5-u^{10})-(u-v)^{12},\]
which also satisfies $P\left(R(q),R\left(q^{11}\right)\right)=0$
\cite{LJRg}. Even if rational parameterizations existed for $R(q)$
and $R\left(q^{11}\right)$, substituting such formulas into Eq.
\eqref{mahler measure for p=11 r(t)} would probably just make the
identity prohibitively complicated.
\subsection{Computationally useful formulas, and a few related hypergeometric transformations}
\label{computational subsection}
While many methods exist for numerically calculating each of the
four Mahler measures $\{\mu(t),n(t),g(t),r(t)\}$, two simple and
efficient methods are directly related to the material discussed so
far.

    The first computational method relies on the $q$-series expansions.  For
example, we can calculate $\mu(\alpha)$ with Eq. \eqref{m(t) q
series}, provided that a value of $q$ exists for which
$M(q)=\alpha$.  Amazingly, the elliptic nome function, defined in
Eq. \eqref{elliptic nome}, furnishes a value of $q$ whenever
$|\alpha|<1$.  Similar inversion formulas exist for all of the
$q$-products in equations \eqref{M(q)} through \eqref{R(q)}. Suppose
that for $j\in\{2,3,4,6\}$
\begin{equation}
q_j(\alpha)=\exp\left(-\frac{\pi}{\sin\left(\frac{\pi}{j}\right)}\frac{{_2F_1}\left(\frac{1}{j},1-\frac{1}{j};1;1-\alpha\right)}
{{_2F_1}\left(\frac{1}{j},1-\frac{1}{j};1;\alpha\right)}\right),
\end{equation}
then we have the following theorem:
\begin{theorem}\label{inversion theorem} With $\alpha$ and $q$ appropriately restricted, the following table gives inversion formulas for
equations \eqref{M(q)} through \eqref{R(q)}:
\begin{center}
    \begin{tabular}
    {|p{1 in}|p{3.5 in}|}
        \hline
        $\alpha$ & $q$ \\
        \hline
        $M(q)$ & $q_2(\alpha)$ \\
        $N(q)$ & $q_3(\alpha)$ \\
        $G(q)$ & $q_2\left(\frac{u(2+u)^3}{(1+2u)^3}\right)$\text{, where $\alpha^3=\frac{u}{2(1+u)^2}$}\\
        $R(q)$ & $q_4\left(\frac{64k\left(1+k-k^2\right)^5}{\left(1+k^2\right)^2\left(\left(1+11k-k^2\right)^2-125k^2\right)^2}\right)$\text{, where $\alpha^5=\frac{k(1-k)^2}{(1+k)^2}$} \\
        \hline
    \end{tabular}
\end{center}
For example: If $|q|<1$ and $\alpha=M(q)$, then $q=q_2(\alpha)$.
\end{theorem}
\begin{proof}  The inversion formulas for $M(q)$ and $G(q)$ follow from
classical eta function identities, and the inversion formula for
$N(q)$ follows from eta function identities in Ramanujan's theory of
signature three.

    The inversion formula for $R(q)$ seems to be new, so we will
prove it.  Let us suppose that $\alpha=R(q)$ and
$k=R(q)R^2\left(q^2\right)$, where $q$ is fixed.  A formula of
Ramanujan \cite{BA} shows that $\alpha^5=\frac{k(1-k)^2}{(1+k)^2}$,
which establishes the second part of the formula.  Now suppose that
$q=q_2(\alpha_2)$, where $\alpha_2=M(q)$. A classical identity shows
that
\begin{equation*}
q\left(-q;q\right)^{24}_{\infty}=\frac{\alpha_2}{16(1-\alpha_2)^2},
\end{equation*}
and comparing this to Ramanujan's identity
\begin{equation*}
q\left(-q;q\right)^{24}_{\infty}=\left(\frac{k}{1-k^2}\right)\left(\frac{1+k-k^2}{1-4k-k^2}\right)^{5},
\end{equation*}
we deduce that
\begin{equation}\label{alpha2 in terms of k}
\frac{\alpha_2}{(1-\alpha_2)^2}=16\left(\frac{k}{1-k^2}\right)\left(\frac{1+k-k^2}{1-4k-k^2}\right)^{5}.
\end{equation}
Now recall that the theory of the signature $4$ elliptic nome shows
that
\begin{equation*}
q=q_2(\alpha_2)=q_4\left(\frac{4\alpha_2}{\left(1+\alpha_2\right)^2}\right)=q_4\left(\frac{4\alpha_2/(1-\alpha_2)^2}{1+4\alpha_2/(1-\alpha_2)^2}\right).
\end{equation*}
Substituting Eq. \eqref{alpha2 in terms of k} into this final result
yields
\begin{equation*}
q=q_4\left(\frac{64k\left(1+k-k^2\right)^5}{\left(1+k^2\right)^2\left(\left(1+11k-k^2\right)^2-125k^2\right)^2}\right),
\end{equation*}
which completes the proof.$\blacksquare$
\end{proof}
The second method for calculating the four Mahler measures,
$\{\mu(t), n(t), g(t), r(t)\}$, depends on reformulating them in
terms of hypergeometric functions. For example, Rodriguez-Villegas
proved \cite{RV} the formula
\begin{equation*}
\mu(t)=-\frac{1}{2}\Re\left[\log(t/16)+\int_{0}^{t}\frac{{_2F_1}\left(\frac{1}{2},\frac{1}{2};1;u\right)-1}{u}\d
u\right].
\end{equation*}
Translated into the language of generalized hypergeometric
functions, this becomes
\begin{equation}\label{m(t) as a 4F3}
\mu(t)=-\Re\left[\frac{t}{8}{_4F_3}\left(\substack{\frac{3}{2},\frac{3}{2},1,1\\2,2,2};t\right)+\frac{1}{2}\log(t/16)\right].
\end{equation}
He also proved a formula for $n(t)$ which is equivalent to
\begin{equation}\label{n(t) as a 4F3}
n(t)=-\Re\left[\frac{2t}{27}{_4F_3}\left(\substack{\frac{4}{3},\frac{5}{3},1,1\\2,2,2};t\right)+\frac{1}{3}\log(t/27)\right].
\end{equation}
Formulas like Eq. \eqref{m(t) as a 4F3} and Eq. \eqref{n(t) as a
4F3} hold some obvious appeal. From a computational perspective they
are useful because most mathematics programs have routines for
calculating generalized hypergeometric functions. For example, when
$|t|<1$ the Taylor series for the ${_4F_3}$ function easily gives
better numerical accuracy than the Mahler measure integrals.
Combining Eq. \eqref{n(t) as a 4F3} with Eq. \eqref{functional
equation g in terms of n} also yields a useful formula for
calculating $g(t)$ whenever $|t|$ is sufficiently small:
\begin{equation}\label{g(t) as a 4F3}
\begin{split}
g(t)=&-\Re\left[
\frac{2t}{(1+4t)^3}{_4F_3}\left(\substack{\frac{4}{3},\frac{5}{3},1,1\\2,2,2};\frac{27t}{(1+4t)^3}\right)
+\frac{8t^2}{(1-2t)^3}{_4F_3}\left(\substack{\frac{4}{3},\frac{5}{3},1,1\\2,2,2};\frac{27t^2}{(1-2t)^3}\right)\right.\\
&\quad\qquad\left.+\log\left(\frac{t^3}{(1+4t)(1-2t)^4}\right)\right].
\end{split}
\end{equation}
So far we have been unable to to find a similar expression for $r(t)$. \\

\textbf{Open Problem 2:} Express $r(t)$ in terms of generalized
hypergeometric functions.\\

Besides their computational importance, identities like Eq.
\eqref{m(t) as a 4F3} allow for a reformulation of Boyd's
conjectures in the language of hypergeometric functions.  For
example, the conjecture
\begin{equation*}
\m\left(1+x+\frac{1}{x}+y+\frac{1}{y}\right)\stackrel{?}{=}L'\left(E,0\right),
\end{equation*}
where $E$ is an elliptic curve with conductor $15$, becomes
\begin{equation*}\label{deninger rephrased}
L'(E,0)\stackrel{?}{=}-2\Re\left[{_4F_3}\left(\substack{\frac{3}{2},\frac{3}{2},1,1\\2,2,2};16\right)\right].
\end{equation*}
A proof of this identity would certainly represent an important
addition to the vast literature concerning transformations and
evaluations of generalized hypergeometric functions.

In the remainder of this section we will apply our results to deduce
a few interesting hypergeometric transformations.  For example,
differentiating Eq. \eqref{g(t) as a 4F3} leads to an interesting
corollary:
\begin{corollary} For $|t|$ sufficiently small
\begin{equation}\label{omega closed form}
\begin{split}
\omega(t):=\sum_{n=0}^{\infty}t^n\sum_{k=0}^{n}{n\choose
k}^3=\frac{1}{1-2t}{_2F_1}\left(\frac{1}{3},\frac{2}{3};1;\frac{27t^2}{(1-2t)^3}\right),
\end{split}
\end{equation}
furthermore
\begin{equation}\label{omega functional equation}
\omega\left(\frac{p}{2(1+p)^2}\right)=(1+p)\omega\left(\frac{p^2}{4(1+p)}\right),
\end{equation}
whenever $|p|$ is sufficiently small.
\end{corollary}
\begin{proof} We can prove Eq. \eqref{omega closed form} by differentiating each side of Eq. \eqref{g(t) as a 4F3}, and then by appealing to Stienstra's
formulas \cite{St}.  A second possible proof follows from showing
that both sides of  Eq. \eqref{omega closed form} satisfy the same
differential equation.

The shortest proof of Eq. \eqref{omega functional equation} follows
from a formula due to Zagier \cite{St}:
\begin{equation*}
\omega\left(G^3(q)\right)=\prod_{n=0}^{\infty}\frac{\left(1-q^{2n}\right)\left(1-q^{3n}\right)^6}
{\left(1-q^n\right)^2\left(1-q^{6n}\right)^3}.
\end{equation*}
First use Zagier's identity to verify that
$G^2(q)\omega\left(G^3(q)\right)=G\left(q^2\right)\omega\left(G^3\left(q^2\right)\right)$,
and then apply the parameterizations for $G^3(q)$ and
$G^3\left(q^2\right)$ from Theorem \ref{functional equation n and g
theorem}.$\blacksquare$
\end{proof}

We will also make a few remarks about the derivative of $r(t)$.
Stienstra has shown that
\begin{equation}\label{psi(t) integral involving phi(t)}
r(t)=-\Re\left[\log(t)+\int_{0}^{t}\frac{\phi(u)-1}{u}\d u\right],
\end{equation}
where $\phi(t)$ is defined by
\begin{equation}\label{phi(t) definition}
\phi(t)=\sum_{n=0}^{\infty}t^n\sum_{k=0}^{n}{n\choose
k}^2{n+k\choose k}.
\end{equation}
Even though we have not discovered a formula for $r(t)$ involving
hypergeometric functions, we can still express $\phi(t)$ in terms of
the hypergeometric function.
\begin{theorem}\label{phi(t) theorem} Let $\phi(t)$ be defined by Eq. \eqref{phi(t) definition},
then for $|k|$ sufficiently small:
\begin{align}
\begin{split}
\phi\left(k\left(\frac{1-k}{1+k}\right)^2\right)=&\frac{(1+k)^2}
{\sqrt{\left(1+k^2\right)\left(\left(1-k-k^2\right)^2-5k^2\right)}}\\
&\times{_2F_1}\left(\frac{1}{4},\frac{3}{4};1;\frac{64k^5\left(1+k-k^2\right)}{\left(1+k^2\right)^2\left(\left(1-k-k^2\right)^2-5k^2\right)^2}\right),
\end{split}\label{phi(t) closed form}\\
\begin{split}
\phi\left(k^2\left(\frac{1+k}{1-k}\right)\right)=&\frac{(1-k)}
{\sqrt{\left(1+k^2\right)\left(\left(1+11k-k^2\right)^2-125k^2\right)}}\\
&\times{_2F_1}\left(\frac{1}{4},\frac{3}{4};1;\frac{64k\left(1+k-k^2\right)^5}{\left(1+k^2\right)^2\left(\left(1+11k-k^2\right)^2-125k^2\right)^2}\right).
\end{split}\label{phi(t) second closed form}
\end{align}

Furthermore, $\phi(t)$ satisfies the functional equation:
\begin{equation}\label{phi(t) cubic functional equation}
\phi\left(k^2\left(\frac{1+k}{1-k}\right)\right)=\frac{1-k}{(1+k)^2}\phi\left(k\left(\frac{1-k}{1+k}\right)^2\right).
\end{equation}
\end{theorem}
\begin{proof}We will prove Eq. \eqref{phi(t) cubic
functional equation} first.  A result of Verrill \cite{Ve} shows
that
\begin{equation}\label{verrill's equation}
\phi^2\left(R^5(q)\right)=\frac{q}{R^5(q)}\frac{\left(q^5;q^5\right)^5_{\infty}}{\left(q;q\right)_{\infty}}.
\end{equation}
Combining Eq. \eqref{verrill's equation} with the trivial formula
$(q^2,q^2)_\infty=(q;q)_{\infty}(-q;q)_{\infty}$, we have
\begin{equation}\label{phi(t) q-ratio}
\frac{\phi^2\left(R^5(q)\right)}{\phi^2\left(R^5(q^2)\right)}=
\frac{R^5\left(q^2\right)}{R^5(q)}\frac{\left\{q^{1/24}\left(-q;q\right)_{\infty}\right\}}{\left\{q^{5/24}\left(-q^{5};q^{5}\right)_{\infty}\right\}^5}.
\end{equation}
We will apply four of Ramanujan's formulas to finish the proof. If
$k=R(q)R^2(q^2)$, then for $|q|$ sufficiently small \cite{BA}:
\begin{align}
R^5(q)=&k\left(\frac{1-k}{1+k}\right)^2,\label{R^5(q)
parameterized}\\
R^5(q^2)=&k^2\left(\frac{1+k}{1-k}\right)\label{R^5(q^2) parameterized},\\
q^{1/24}\left(-q;q\right)_{\infty}=&\left(\frac{k}{1-k^2}\right)^{1/24}\left(\frac{1+k-k^2}{1-4k-k^2}\right)^{5/24},\label{R(q)
(-q,q)infinity parameterization}\\
q^{5/24}\left(-q^5;q^5\right)_{\infty}=&\left(\frac{k}{1-k^2}\right)^{5/24}\left(\frac{1+k-k^2}{1-4k-k^2}\right)^{1/24}.\label{R(q)
(-q^5,q^5)infinity parameterization}
\end{align}
Eq. \eqref{phi(t) cubic functional equation} follows immediately
from substituting these parametric formulas into Eq. \eqref{phi(t)
q-ratio}.

    Next we will prove Eq. \eqref{phi(t) closed form}.
Combining Eq. \eqref{R^5(q) parameterized} with Entry 3.2.15 in
\cite{BA}, we easily obtain
\begin{equation}\label{phi multiplier}
q^{5/24}\left(q^5;q^5\right)_{\infty}
=\left\{\frac{k(1-k^2)^2}{\left(1+k-k^2\right)\left(1-4k-k^2\right)^2}\right\}^{1/6}q^{1/24}\left(q;q\right)_{\infty}.
\end{equation}
Now we will evaluate the eta product
$q^{1/24}\left(q;q\right)_{\infty}$. First recall that if
$q=q_4(z)$, then
\[q^{1/24}(q;q)_{\infty}=2^{-1/4}z^{1/24}(1-z)^{1/12}\sqrt{{_2F_1}\left(\frac{1}{4},\frac{3}{4};1;z\right)}.\]
In Theorem \ref{inversion theorem} we showed that if
$k=R(q)R^2\left(q^2\right)$ then
$q=q_4\left(\frac{64k\left(1+k-k^2\right)^5}{\left(1+k^2\right)^2\left(\left(1+11k-k^2\right)^2-125k^2\right)^2}\right)$,
hence it follows that
\begin{equation}\label{phi eta^5  eval}
\begin{split}
q^{1/24}\left(q;q\right)_{\infty}=
&\left(\frac{k\left(1-k^2\right)^{2}
\left(1+k-k^2\right)^5\left(1-4k-k^2\right)^{10}}{\left(1+k^2\right)^{6}\left(\left(1+11 k-k^2\right)^2-125k^2\right)^6}\right)^{1/24}\\
&\times\sqrt{{_2F_1}\left(\frac{1}{4},\frac{3}{4};1;\frac{64k\left(1+k-k^2\right)^5}{\left(1+k^2\right)^2\left(\left(1+11k-k^2\right)^2-125k^2\right)^2}\right)}
\end{split}
\end{equation}
Substituting Eq. \eqref{phi eta^5  eval}, Eq. \eqref{phi
multiplier}, and Eq. \eqref{R^5(q)
    parameterized} into Eq. \eqref{verrill's equation}
completes the proof of Eq. \eqref{phi(t) closed form}.  The proof of
Eq. \eqref{phi(t) second closed form} also follows from an extremely
similar argument.
 $\blacksquare$
\end{proof}
We will conclude this section by recording a few formulas which do
not appear in \cite{BA}, but which were probably known to Ramanujan.
We will point out that  Maier obtained several results along these
lines in \cite{Ma}.  Notice that the functional equation for
$\phi(t)$ (after substituting $z=k/(1-k^2)$) implies a new
hypergeometric transformation:
\begin{equation}
\begin{split}
\sqrt{\frac{(1+11z)^2-125z^2}{(1-z)^2-5
z^2}}&{_2F_1}\left(\frac{1}{4},\frac{3}{4};1,\frac{64z^5\left(1+z\right)}{\left(1+4z^2\right)\left(\left(1-z\right)^2-5z^2\right)^2}\right)\\
&={_2F_1}\left(\frac{1}{4},\frac{3}{4};1;\frac{64z\left(1+z\right)^5}{\left(1+4z^2\right)\left(\left(1+11z\right)^2-125z^2\right)^2}\right)
\end{split}
\end{equation}
Perhaps not surprisingly, we can also use the arguments in this
section to deduce that
\begin{equation}
\begin{split}
q_4^5\left(\frac{64z\left(1+z\right)^5}{\left(1+4z^2\right)\left(\left(1+11z\right)^2-125z^2\right)^2}\right)=q_4\left(\frac{64z^5\left(1+z\right)}{\left(1+4z^2\right)\left(\left(1-z\right)^2-5z^2\right)^2}\right),
\end{split}
\end{equation}
which implies a rational parametrization for the fifth-degree
modular equation in Ramanujan's theory of signature $4$.

%
%

\section{A regulator explanation} \label{regul}
Now  we will reinterpret our identities in terms of the regulators
of elliptic curves.  The elliptic curves in question are defined
by the zero varieties of the polynomials whose Mahler measure we
studied. First we will explain the relationship between Mahler
measures and regulators. Then we will use regulators to deduce
formulas involving Kronecker-Eisenstein series, including
equations \eqref{m(t) q series}, \eqref{n(t) q series}, \eqref{g(t) q series}, and \eqref{r(t) q series}.

We will follow some of the ideas of Rodriguez-Villegas \cite{RV2}.

\subsection{The elliptic regulator}
Let $F$ be a field. By Matsumoto's Theorem, $K_2(F)$ is generated
by the symbols $\{a,b\}$ for $a,b\in F$, which satisfy the
bilinearity relations $\{a_1a_2,b\} = \{a_1,b\}\{a_2,b\}$ and
$\{a,b_1b_2\}=\{a,b_1\}\{a,b_2\}$, and the Steinberg relation
$\{a,1-a\}=1$ for all $a \not =0$.

Recall that for a field $F$, with discrete valuation $v$, and
maximal ideal $\mathcal{M}$, the tame symbol is given by
\[ (x,y)_v \equiv (-1)^{v(x)v(y)} \frac{x^{v(y)}}{y^{v(x)}} \,\, \mathrm{mod} \,\, \mathcal{M}\]
(see \cite{RV}). Note that this symbol is trivial if
$v(x)=v(y)=0$. In the case when $F=\mathbb{Q}(E)$ (from now on $E$
denotes an elliptic curve), a valuation is determined by the order
of the rational functions at each point $S\in
E(\bar{\mathbb{Q}})$. We will denote the valuation determined by a
point $S\in E(\bar{\mathbb{Q}})$ by $v_S$.

The tame symbol is then a map $K_2(\mathbb{Q}(E)) \rightarrow \mathbb{Q}(S)^*$.

We have
\[ 0 \rightarrow K_2(E) \otimes \mathbb{Q} \rightarrow
K_2(\mathbb{Q}(E))\otimes\mathbb{Q} \rightarrow \coprod_{S\in
E(\bar{\mathbb{Q}})}\mathbb{Q}Q(S)^*\times \mathbb{Q}, \] where
the last arrow corresponds to the coproduct of tame symbols.

Hence an element $\{x,y\} \in K_2(\mathbb{Q}(E))\otimes
\mathbb{Q}$ can be seen as an element in $K_2(E) \otimes
\mathbb{Q}$ whenever $(x,y)_{v_S}=1$ for all $S \in
E(\bar{\mathbb{Q}})$. All of the families considered in this paper
are tempered according to \cite{RV}, and therefore they satisfy
the triviality of tame symbols.

The regulator map (defined by Beilinson after the work of Bloch) may be defined by
\[\r :K_2(E) \rightarrow H^1(E,\mathbb{R})\]
\[\{x,y\}\rightarrow \left\{ \gamma \rightarrow \int_\gamma \eta(x,y)\right \}\]
for $\gamma \in H_1(E,\mathbb{Z})$, and
\[\eta(x,y) := \log |x| \d \arg y - \log|y| \d \arg x. \]
Here we think of $H^1(E,\mathbb{R})$ as the dual of
$H_1(E,\mathbb{Z})$. The regulator is well defined because
$\eta(x,1-x) = \d D(x)$, where
\[D(z)= \Im(\Li_2(z))+\arg(1-z) \log|z|\]
is the Bloch-Wigner dilogarithm.

In terms of the general formulation of Beilinson's conjectures
this definition is not completely correct.  One needs to go a step
further and consider $K_2(\mathcal{E})$, where $\mathcal{E}$ is a
N\'eron model of $E$ over $\mathbb{Z}$. In particular,
$K_2(\mathcal{E})$ is a subgroup of $K_2(E)$. It seems \cite{RV}
that a power of $\{x,y\}$ always lies in $K_2(\mathcal{E})$.

Assume that $E$ is defined over $\mathbb{R}$. Because of the way
that complex conjugation acts on $\eta$, the regulator map is
trivial for the classes in $H_1(E,\mathbb{Z})^+$.  In particular,
these cycles remain invariant under complex conjugation. Therefore
it suffices to consider the regulator as a function on
$H_1(E,\mathbb{Z})^-$.

We write $E(\mathbb{C}) \cong
 \mathbb{C}/\mathbb{Z}+\tau\mathbb{Z}$, where $\tau$ is in the upper half-plane. Then
 $\mathbb{C}/\mathbb{Z}+\tau\mathbb{Z} \cong \mathbb{C}^*/q^\mathbb{Z}$, where
 $z$ mod $\Lambda= \mathbb{Z}+\tau \mathbb{Z}$ is identified with $\e^{2 \ii\pi
 z}$.
Bloch \cite{B} defines the regulator function  in terms of a
Kronecker-Eisenstein series
\begin{equation}
R_\tau \left(\e^{2\pi\ii (a + b\tau)} \right) =
\frac{y_\tau^2}{\pi} \sum_{m,n \in \mathbb{Z}}'
\frac{\e^{2\pi\ii(bn-am)}}{(m\tau+n)^2(m\bar{\tau}+n)},
\end{equation}
where $y_\tau$ is the imaginary part of $\tau$.

Let $J(z) = \log|z| \log|1-z|$, and let
\[D(x)= \Im(\Li_2(x))+\arg(1-x) \log|x|\]
be the Bloch-Wigner dilogarithm.

Consider the following function on $E(\mathbb{C}) \cong
\mathbb{C}^*/q^\mathbb{Z}$:
\begin{equation}
J_\tau\left(z\right) = \sum_{n=0}^\infty J(zq^n) -
\sum_{n=1}^\infty J(z^{-1}q^n) + \frac{1}{3} \log^2|q|
B_3\left(\frac{\log|z|}{\log|q|} \right),
\end{equation}
where $B_3(x)= x^3-\frac{3}{2}x^2 +\frac{1}{2}x$ is the third
Bernoulli polynomial.  If we recall that the elliptic dilogarithm
is defined by
\begin{equation}
D_\tau(z) := \sum_{n \in \mathbb{Z}} D(z q^n),
\end{equation}
then the regulator function (see \cite{B}) is given by
\begin{equation}
R_\tau = D_\tau -\ii J_\tau.
\end{equation}

By linearity, $R_\tau$ extends to divisors with support in
$E(\mathbb{C})$. Let $x$ and $y$ be non-constant functions on $E$
with divisors
\[(x) = \sum m_i (a_i), \qquad (y) = \sum n_j (b_j).\]
Following \cite{B}, and the notation in \cite{RV}, we recall the
diamond operation $\mathbb{C}(E)^* \otimes  \mathbb{C}(E)^*
\rightarrow \mathbb{Z}[E(\mathbb{C})]^-$
\[ (x)\diamond (y) = \sum m_i n_j (a_i-b_j).\]
Here $\mathbb{Z}[E(\mathbb{C})]^-$ means that $[-P]\sim -[P]$.

Because $R_\tau$ is an odd function, we obtain a
 map \[ \mathbb{Z}[E(\mathbb{C})]^- \rightarrow \mathbb{R}.\]

\begin{theorem} (Beilinson \cite{Be}) $E/\mathbb{R}$ elliptic curve, $x, y $ non-constant
 functions in $\mathbb{C}(E)$, $\omega \in \Omega^1$
\[\int_{E(\mathbb{C})}
\bar{\omega} \wedge \eta(x,y) =  \Omega_0 R_\tau((x)\diamond(y))\]
 where $\Omega_0$ is the real period.
 \end{theorem}
Although a more general version of Beilinson's Theorem exists for
elliptic curves defined over the complex numbers, the above
version has a simpler formulation.

 \begin{corollary} (after an idea of Deninger) If
$x$ and $y$ are non-constant functions in $\mathbb{C}(E)$ with
trivial tame symbols, then
\[
-\int_\gamma \eta(x,y)= \Im\left(\frac{\Omega}{ y_\tau
\Omega_0}
R_\tau\left( (x)\diamond (y) \right) \right)\]
where $\Omega=\int_\gamma \omega$.
\end{corollary}

\begin{proof}  Notice that $\ii \eta(x,y) $ is an element of the two-dimensional vector
space $H^2_{\mathcal{D}}(E(\mathbb{C}),\mathbb{R}(2))$ generated
by $\omega$ and $\bar{\omega}$. Then we may write
\[\ii \eta(x,y) = \alpha [\omega]+\beta [\bar{\omega}],\] from which we obtain \[
\int_\gamma \ii\eta(x,y) = \alpha
\Omega + \beta \overline{\Omega}.\] On the other hand, we have
\[\int_{E(\mathbb{C})} \ii \eta(x,y)\wedge
\bar{\omega} =\alpha  \int_{E(\mathbb{C})} \omega \wedge
\bar{\omega} = \alpha  \ii 2 \Omega_0^2 y_\tau, \] and \[
\int_{E(\mathbb{C})} \ii \eta(x,y)\wedge \omega =- \beta  \ii 2
\Omega_0^2 y_\tau. \] By Beilinson's Theorem
\[\int_\gamma \ii \eta(x,y)=-\frac{R_\tau((x)\diamond (y)) \Omega}{ 2 \Omega_0
y_\tau} + \frac{\overline{R_\tau((x)\diamond (y))}
\overline{\Omega}}{ 2 \Omega_0 y_\tau},\] and the statement
follows.$\blacksquare$
\end{proof}



\subsection{Regulators and Mahler measure}
From now on, we will set $k=\frac{4}{\sqrt{t}}$ in the first family \eqref{eq:m}.

Rodriguez-Villegas \cite{RV} proved that if $P_k(x,y)=k + x +
\frac{1}{x} + y + \frac{1}{y}$ does not intersect the torus
$\mathbb{T}^2$, then
\begin{equation}\label{eq:regmm}
m(k)\sim_{\mathbb{Z}} \frac{1}{2\pi} \r(\{x,y\})(\gamma).
\end{equation}
Here the $\sim_{\mathbb{Z}}$ stands for "up to an integer number",
and $\gamma$ is a closed path that avoids the poles and zeros of
$x$ and $y$.  In particular, $\gamma$ generates the subgroup
$H_1(E,\mathbb{Z})^-$ of $H_1(E,\mathbb{Z})$ where conjugation
acts by $-1$.

We would like to use this property, however we need to exercise
caution.  In particular, $P_k(x,y)$ intersects the torus whenever
$|k|\leq 4$ and $k\in \mathbb{R}$. Let us recall the idea behind
the proof of Eq. \eqref{eq:regmm} for the special case of
$P_k(x,y)$. Writing
\[yP_k(x,y) = (y-y_{(1)}(x))(y-y_{(2)}(x)),\]
we have
\[m(k) = \m(y P_k(x,y)) = \frac{1}{2\pi\ii} \int_{\mathbb{T}^1} (\log^+|y_{(1)}(x)|+\log^+|y_{(2)}(x)|) \frac{\d x}{x}. \]
This last equality follows from applying Jensen's formula with
respect to the variable $y$. When the polynomial does not
intersect the torus, we may omit the ``$+$" sign on the logarithm
since each $y_{(i)}(x)$ is always inside or outside the unit
circle. Indeed, there is always a branch inside the unit circle
and a branch outside. It follows that
\begin{equation}\label{eq:mm}
 m(k)= \frac{1}{2\pi\ii} \int_{\mathbb{T}^1} \log |y| \frac{\d x}{x} = -\frac{1}{2\pi} \int_{\mathbb{T}^1} \eta(x,y),
 \end{equation}
where $\mathbb{T}^1$ is interpreted as a cycle in the homology of
the elliptic curve defined by $P_k(x,y)=0$, namely
$H_1(E,\mathbb{Z})$.

If $k\in[-4,4]$, then we may also assume that $k>0$ since this
particular Mahler measure does not depend on the sign of $k$. The
equation
\[ k+x+\frac{1}{x}+y+\frac{1}{y} =0 \]
certainly has solutions when $(x,y)\in \mathbb{T}^2$. However, for
$|x|=1$ and $k$ real, the number $k+x+\frac{1}{x}$ is real, and
therefore $y+\frac{1}{y}$ must be real. This forces two
possibilities: either $y$ is real or $|y|=1$. Let $x=\e^{\ii
\theta}$, then for $-\pi\leq \theta\leq \pi$ we have
\begin{equation}\label{eq:Y}
 -k-2 \cos \theta = y +\frac{1}{y}.
 \end{equation}
The limiting case occurs when $|k+2\cos \theta| = 2$. Since we
have assumed that $k$ is positive, this condition becomes $k+2\cos
\theta=2$, which implies that $y=-1$. When $k+2\cos \theta
> 2$ one solution for $y$, say, $y_{(1)}$, becomes a negative
number less than -1, thus $\left|y_{(1)}\right| >1$ (the other
solution $y_{(2)}$ is such that $\left|y_{(2)}\right|<1$). When
$k+2\cos \theta < 2$, $y_i$ lies inside the unit circle and never
reaches 1. What is important is that $\left|y_{(1)}\right|\geq 1$
and $\left|y_{(2)}\right|\leq 1$, so we can still write Eq.
\eqref{eq:mm} even if there is a nontrivial intersection with the
torus.

\subsection{Functional identities for the regulator}

First recall a result by Bloch \cite{B} which studies the
modularity of $R_\tau$:

\begin{proposition} \label{modular prop of R}
Let $\left( \begin{array}{cc} \alpha & \beta \\ \gamma & \delta
\end{array} \right) \in SL_2(\mathbb{Z})$, and let $\tau' =
\frac{\alpha \tau + \beta}{\gamma \tau + \delta} $. If we let
\[\left( \begin{array}{c} b' \\ a' \end{array} \right) = \left( \begin{array}{cc} \delta & -\gamma \\ -\beta & \alpha \end{array} \right)\left( \begin{array}{c} b \\ a \end{array} \right), \]
then:
\[R_{\tau'}\left(\e^{2\pi\ii (a'+b'\tau')}\right)=\frac{1}{\gamma \bar{\tau} +\delta} R_{\tau}\left(\e^{2\pi\ii(a+b\tau)}\right).\]
\end{proposition}\allowdisplaybreaks{

We will need to use some functional equations for $J_\tau$. First
recall the following trivial property for $J(z)$:
\begin{equation}\label{J-equation}
J(z)=p\sum_{x^p=z}J(x).
\end{equation}

\begin{proposition}
 Let $p$ be an odd prime, let $q=\e^{2\pi \ii \tau}$,
 and let $q_j=\e^{\frac{2\pi \ii (\tau+j)}{p}}$ for $j\in\{0,1,\dots, p-1\}$.
 Suppose that $(N,k)=1$, and $p\equiv \pm1$ or $0\, \left(\mathrm{mod}\, N\right)$. Then
 \begin{equation} \label{eq:JNk}
 (1+\chi_{-N}(p)p^2) J_{N\tau}(q^k) = \sum_{j=0}^{p-1}p J_{\frac{N(\tau+j)}{p}}(q_j^k) +\chi_{-N}(p)
 J_{Np\tau}(q^{pk}),
 \end{equation}
and for any $z$ we have
  \begin{equation}
  (\chi_{-N}(p)+ p^2) J_{N\tau}(z) = \sum_{j=0}^{p-1}p J_{\frac{N(\tau+j)}{p}}(z) +\chi_{-N}(p)
  J_{Np\tau}(z).
 \end{equation}
%
\end{proposition}
\begin{proof}\allowdisplaybreaks{
First notice that
\begin{equation*}
\begin{split}
\sum_{j=0}^{p-1} J_{\frac{N(\tau+j)}{p}}(q_j^k) =&
\sum_{n=0}^\infty \sum_{j=0}^{p-1} J\left(q_j^{Nn+k}\right)\\
& - \sum_{n=1}^\infty \sum_{j=0}^{p-1} J\left(q_j^{Nn-k}\right) +
\frac{4\pi^2y_\tau^2 N^2}{3p} B_3\left(\frac{k}{N}\right).
\end{split}
\end{equation*}
 By Eq. \eqref{J-equation} this becomes
\begin{align*}
=&\sum_{\substack{n=0\\p\nmid Nn+k} }^\infty \frac{1}{p}
J\left(q^{Nn+k}\right) - \sum_{\substack{n=1\\ p \nmid
Nn-k}}^\infty
\frac{1}{p} J\left(q^{Nn-k}\right)\\
&+\sum_{\substack{n=0\\ p | Nn+k} }^\infty p
J\left(q^{\frac{Nn+k}{p}}\right) - \sum_{\substack{n=1\\ p |
Nn-k}}^\infty p J\left(q^{\frac{Nn-k}{p}}\right)
+\frac{4\pi^2y_\tau^2 N^2}{3p}
B_3\left(\frac{k}{N}\right)\\
 =& \sum_{n=0}^\infty \frac{1}{p} J\left(q^{Nn+k}\right) - \sum_{n=1 }^\infty \frac{1}{p}
 J\left(q^{Nn-k}\right)\\
  &-\sum_{\substack{n=0\\ p | Nn+k} }^\infty \frac{1}{p}
J\left(q^{Nn+k}\right) + \sum_{\substack{n=1\\ p | Nn-k}}^\infty
\frac{1}{p} J\left(q^{Nn-k}\right)\\
&+\sum_{\substack{n=0\\ p | Nn+k} }^\infty p
J\left(q^{\frac{Nn+k}{p}}\right) - \sum_{\substack{n=1\\ p |
Nn-k}}^\infty p J\left(q^{\frac{Nn-k}{p}}\right)
+\frac{4\pi^2y_\tau^2 N^2}{3p} B_3\left(\frac{k}{N}\right).
\end{align*}
Rearranging, we find that
\begin{align*}
 =&\frac{1}{p} J_{{N\tau}}\left(q^k\right) -
\frac{4\pi^2y_\tau^2 N^2}{3p} B_3\left(\frac{k}{N}\right)\\
&- \sum_{\substack{n=0\\ p | Nn+k} }^\infty \frac{1}{p}
J\left((q^p)^{\frac{Nn+k}{p}}\right) + \sum_{\substack{n=1\\ p |
Nn-k}}^\infty \frac{1}{p} J\left((q^p)^{\frac{Nn-k}{p}}\right)\\
&+\sum_{\substack{n=0\\ p | Nn+k} }^\infty p
J\left(q^{\frac{Nn+k}{p}}\right) - \sum_{\substack{n=1\\ p |
Nn-k}}^\infty p J\left(q^{\frac{Nn-k}{p}}\right)
+\frac{4\pi^2y_\tau^2 N^2}{3p}
B_3\left(\frac{k}{N}\right)\\
= & \frac{1}{p} J_{{N\tau}}\left(q^k\right) -
\frac{\chi_{-N}(p)}{p} J_{Np\tau}(q^{pk}) + \chi_{-N}(p)p
J_{N\tau}(q^k),
\end{align*}
which proves the assertion.

The second equality follows in a similar fashion. $\blacksquare$}
 \end{proof}
It is possible to prove analogous identities for $D_\tau$ and
$R_\tau$.

\begin{proposition}
\begin{equation} \label{eq:J2}
 J_{\frac{2\mu+1}{2}}\left( \e^{\pi\ii \mu }\right) = J_{2\mu}\left(\e^{\pi\ii \mu}\right) - J_{2\mu}\left(-\e^{\pi\ii \mu}\right)
 \end{equation}
 \end{proposition}
\begin{proof}
Let $z=\e^{\pi\ii \mu}$, then
\begin{align*}
J_{2\mu}\left(z\right) - J_{2\mu}\left(-z\right)=&J(z)-J(-z)\\
&+\sum_{n=1}^\infty \left(J(zq^n) -J(-zq^n)-J(z^{-1}q^n)+J(-z^{-1}q^n)\right)\\
=&\sum_{n=0}^\infty \left(J\left(\e^{\pi\ii \mu (4n+1)}\right)-
J\left(-\e^{\pi\ii \mu (4n+1)} \right)\right.\\
&\qquad\left.- J\left(\e^{\pi\ii \mu (4n+3)} \right) +
J\left(-\e^{\pi\ii \mu (4n+3)} \right)\right).
\end{align*}
On the other hand,
\[ J_{\frac{2\mu+1}{2}}\left(z\right) = \sum_{n=0}^\infty \left(J\left((-1)^{n} \e^{\pi\ii\mu(2n+1)}\right) - J\left((-1)^{n+1} \e^{\pi\ii\mu(2n+1)}\right)\right),\]
which proves the equality.$\blacksquare$
\end{proof}

\subsection{The first family}

First we will write the equation
\[x+\frac{1}{x}+y+\frac{1}{y}+k =0 \]
in Weierstrass form. Consider the rational transformation
\[ X=\frac{k+x+y}{x+y} = -\frac{1}{xy}, \qquad Y=\frac{k(y-x)(k+x+y)}{2(x+y)^2} = \frac{(y-x)\left(1+\frac{1}{xy}\right)}{2xy}, \]
which leads to
\[ Y^2=X\left(X^2+\left(\frac{k^2}{4}-2\right)X+1\right).\]
It is useful to state the inverse transformation:
\[ x=\frac{kX-2Y}{2X(X-1)}, \qquad y=\frac{kX+2Y}{2X(X-1)}.\]

Notice that $E_k$ contains a torsion point of order $4$ over $\mathbb{Q}(k)$,
namely $P=\left(1, \frac{k}{2}\right)$. Indeed, this family is the modular elliptic surface associated to $\Gamma_0(4)$.

We can show that
$2P=(0,0)$, and $3P=\left(1,-\frac{k}{2}\right)$.

Now \[(X)= 2(2P) -2O,\] and
\begin{equation*}
\begin{split}
(x)=& (2(P)+(2P)-3O) - (2(2P)-2O)-((P)+(3P)-2O)\\
   =& (P)-(2P)-(3P)+O,
\end{split}
\end{equation*}
\begin{equation*}
\begin{split}
(y)=& (2(3P)+(2P)-3O) - (2(2P)-2O)-((P)+(3P)-2O)\\
   =&-(P)-(2P)+(3P)+O.
\end{split}
\end{equation*}
Computing the diamond operation between the divisors of $x$ and
$y$ yields
\[(x)\diamond (y) = 4 (P) -4(-P) = 8(P).\]

Now assume that $k \in \mathbb{R}$ and $k > 4$. We will choose an
orientation for the curve and compute the real period. Because $P$
is a point of order 4 and $\int_0^1 \omega$ is real, we may assume
that $P$ corresponds to $\frac{3\Omega_0}{4}$.

%
%

The next step is to understand the cycle $|x|=1$ as an element of
$H_1(E,\mathbb{Z})$.  We would like to compute the value of
$\Omega = \int_\gamma \omega$. First recall that
\[ \omega = \frac{\d X}{2Y} = \frac{\d x}{x(y-y^{-1})}.\]
In the case when $k>4$, consider conjugation of $\omega$. This
sends $x \rightarrow x^{-1}$, and $\frac{\d x}{x} \rightarrow
-\frac{\d x}{x}$. There is no intersection with the torus, so $y$
remains invariant.  Therefore we conclude that $\Omega $ is the
complex period, and $\frac{\Omega}{\Omega_0} = \tau$, where $\tau$
is purely imaginary.

%
%
%
%
%
Therefore for $k$ real and $|k| >4$
\[m(k) =  \frac{4}{\pi} \Im \left(\frac{\tau}{y_\tau}R_\tau(-\ii) \right).\]
%
%
 Now take $\left(\begin{array}{cc} 0& -1 \\ 1 & 0 \end{array}
\right) \in SL_2(\mathbb{Z})$. By Proposition \ref{modular prop of
R}
\[ R_\tau(-\ii) =R_\tau\left(\e^{-\frac{2\pi
\ii}{4}}\right)=\bar{\tau}
R_{-\frac{1}{\tau}}\left(\e^{-\frac{2\pi \ii}{4\tau}}\right),\]
therefore
 \[  m(k) = - \frac{4 |\tau|^2}{\pi y_\tau}
J_{-\frac{1}{\tau}}\left(\e^{-\frac{2\pi \ii}{4\tau}}\right). \]
If we let $\mu = -\frac{1}{4 \tau}$, then for $k\in \mathbb{R}$ we
obtain
\begin{align*}
 m(k) =& - \frac{1}{\pi y_\mu} J_{4 \mu}\left(\e^{2\pi \ii \mu
}\right) =\Im \left( \frac{1}{\pi y_\mu} R_{4 \mu}\left(\e^{2\pi
\ii \mu }\right) \right)\\
 =&  \Re \left( \frac{16 y_\mu
}{\pi^2} \sum_{m,n}' \frac{\chi_{-4}\left(m\right)}{(m+4\mu
n)^2(m+4\bar{\mu}n)}\right),
\end{align*}
thus recovering a result of Rodriguez-Villegas.
%
%
We can extend this result to all $k \in \mathbb{C}$, by arguing
that both $m(k)$ and $ - \frac{1}{\pi y_\mu} J_{4
\mu}\left(\e^{2\pi \ii \mu }\right)$ are the real parts of
holomorphic functions that coincide at infinitely many points (see
\cite{RV0}).

%

Now we will show how to deduce equations \eqref{eq:first} and \eqref{eq:ko}.
Applying Eq. \eqref{eq:JNk} with $N=4$, $k=1$, and $p=2$, we have
\[  J_{4\mu}(q) = 2 J_{2\mu}(q_0) + 2 J_{2(\mu+1)}(q_1),  \]
which  translates into
\[  \frac{1}{y_{4\mu}} J_{4\mu}\left(\e^{2\pi\ii \mu}\right) = \frac{1}{y_{2\mu}}
J_{2\mu}\left(\e^{\pi\ii \mu}\right) + \frac{1}{y_{2\mu}}
J_{2\mu}\left(-\e^{\pi\ii \mu}\right) .  \] This is the content of
Eq. \eqref{eq:first}. Setting
$\tau=-\frac{1}{2\mu}$, we may also write
\begin{equation}\label{eq:1}
 D_\frac{\tau}{2}(-\ii) =D_{\tau}(-\ii) + D_{\tau}(-\ii \e^{\pi\ii \tau}).
 \end{equation}

Next we will use Eq. \eqref{eq:J2}:
\[ J_{\frac{2\mu+1}{2}}\left( \e^{\pi\ii \mu }\right) = J_{2\mu}\left(\e^{\pi\ii
\mu}\right) - J_{2\mu}\left(-\e^{\pi\ii \mu}\right),\] which
translates into
\[ \frac{1}{y_{\frac{2\mu+1}{2}}} J_{\frac{2\mu+1}{2}}\left( \e^{\pi\ii \mu }\right) = \frac{2}{y_{2\mu}}J_{2\mu}\left(\e^{\pi\ii
\mu}\right) - \frac{2}{y_{2\mu}}J_{2\mu}\left(-\e^{\pi\ii
\mu}\right).\]  Setting $\tau=-\frac{1}{2\mu}$, and using
$\left(\begin{array}{cc} 1 &  0\\ -2 & 1 \end{array} \right) \in
SL_2(\mathbb{Z})$  on the left-hand side, we have
\begin{equation}\label{eq:2}
D_\frac{\tau-1}{2}(-\ii) =  D_\tau(-\ii) -
D_\tau\left(-\ii\e^{\pi\ii \tau} \right).
\end{equation}
Combining equations \eqref{eq:1} and \eqref{eq:2}, we see that
\[2 D_\tau(-\ii) = D_\frac{\tau}{2}(-\ii) + D_\frac{\tau-1}{2}(-\ii).\]
This is the content of Eq. \eqref{eq:ko}.

Similarly, we may deduce Eq. \eqref{Hecke eigenvalue equation}
from Eq. \eqref{eq:JNk} when $k=1$, $N=4$, and $p$ is an odd
prime.

\subsection{A direct approach}
It is also possible to prove equations \eqref{eq:ko} and
\eqref{eq:first} directly, without considering the
$\mu$-parametrization or the explicit form of the regulator.

For those formulas, it is easy to explicitly write the isogenies
at the level of the Weierstrass models. By using the well-known
isogeny of degree 2 (see for example \cite{C,S}):
\[\phi: \{E: y^2=x (x^2+ax+b)\} \rightarrow \{\widehat{E}: \hat{y}^2 = \hat{x}(\hat{x}^2 - 2a \hat{x}+(a^2-4b))\}\]
given by
\[(x,y) \rightarrow \left(\frac{y^2}{x^2}, \frac{y(b-x^2)}{x^2}\right)\]
(we require that $a^2-4b \not = 0$), we find

\[ \phi_1: E_{2\left(n+\frac{1}{n}\right)} \rightarrow E_{4 n^2}, \qquad \phi_2: E_{2\left(n+\frac{1}{n}\right)} \rightarrow E_{\frac{4}{n^2}},\]

\[\phi_1 : \left(X,Y \right) \rightarrow \left( \frac{X(n^2X+1)}{X+n^2}, -\frac{n^3Y\left(X^2+2n^2X+1 \right)}{\left(X+n^2\right)^2}\right),\]
\[\phi_2 : \left(X,Y \right) \rightarrow \left( \frac{X(X+n^2)}{n^2X+1}, -\frac{Y\left(n^2X^2+2X+n^2 \right)}{n\left(n^2X+1\right)^2}\right).\]

Let us write $x_1$, $y_1$, $X_1$, $Y_1$ for the rational functions
and $\r_1$ for the regulator in $E_{4 n^2}$, and $x_2$, $y_2$,
$X_2$, $Y_2$, $\r_2$ for the corresponding objects in
$E_{\frac{4}{n^2}}$.

It follows that
\begin{align*}
\pm m\left(4n^2\right)=& \r_1 \left( \left\{ x_1,y_1 \right\}
\right) = \frac{1}{2\pi} \int_{|X_1|=1} \eta(x_1,y_1)\\
 =& \frac{1}{4\pi} \int_{|X|=1}
\eta(x_1\circ\phi_1,y_1\circ\phi_1)\\
 =& \frac{1}{2}\r \left(
\left\{ x_1\circ\phi_1,y_1\circ\phi_1 \right\} \right),
\end{align*}
where the factor of $2$ follows from the degree of the isogeny.
Similarly, we find that
\[\pm m\left(\frac{4}{n^2}\right)= \r_2 \left( \left\{ x_2,y_2 \right\} \right) =\frac{1}{2}\r \left( \left\{ x_2\circ\phi_2,y_2\circ\phi_2 \right\} \right).\]

Now we need to compare the values of
\[  \r \left( \left\{ x_1\circ\phi_1,y_1\circ\phi_1 \right\} \right),\quad \r \left( \left\{ x_2\circ\phi_2,y_2\circ\phi_2 \right\} \right), \quad \mbox{and} \quad  \r \left( \left\{ x,y \right\} \right).\]

Recall that $(x)\diamond (y) = 8(P)$, where $P=(1, \frac{k}{2})$.
When $k=2\left(n+\frac{1}{n}\right)$, we will also consider the
point $Q=\left(-\frac{1}{n^2},0\right)$, which has order $2$ (then
$P+Q = \left(-1,n-\frac{1}{n}\right)$, $2P+Q
=\left(-n^2,0\right)$, etc).

Let $P$ now denote the point in $E_{2\left(n+\frac{1}{n}\right)}$,
and let $P_1$ denote the corresponding point in $E_{4n^2}$. We
have the following table:
\[ \phi_1: \qquad \begin{array}{cccc} 3P, &P+Q & \rightarrow & P_1\\ 2P, &Q & \rightarrow &2P_1\\ P, &3P+Q & \rightarrow &3P_1\\ O_0, &2P+Q & \rightarrow &O_1  \end{array}.\]
Using this table, and  the divisors $(x_1)$ and $(y_1)$ in
$E_{4n^2}$, we can compute $(x_1\circ\phi_1)\diamond
(y_1\circ\phi_1)$.  We find that
\[(x_1\circ\phi_1)\diamond (y_1\circ\phi_1) = -16(P)+16(P+Q),\]
and similarly
\[(x_2\circ\phi_2)\diamond (y_2\circ\phi_2) = -16(P) -16(P+Q).\]

These computations show that
\begin{equation}
\frac{1}{2} \r_0 \left( \left\{ x_1\circ\phi_1,y_1\circ\phi_1
\right\} \right)  + \frac{1}{2} \r_0 \left( \left\{
x_2\circ\phi_2,y_2\circ\phi_2 \right\} \right) = 2 \r_0 \left(
\left\{ x_0,y_0 \right\} \right),
\end{equation}
and therefore
\begin{equation}
\r_1 \left( \left\{ x_1,y_1\right\} \right)  +  \r_2 \left( \left\{
x_2,y_2 \right\} \right) = 2 \r_0 \left( \left\{ x_0,y_0 \right\}
\right).
\end{equation}
We can conclude the proof of Eq. \eqref{eq:ko} by
inspecting signs.

To prove Eq. \eqref{eq:first}, it
is necessary to use the isomorphism $\phi$ from Eq.
\eqref{eq:iso}.

\subsection{Relations among $m(2)$, $m(8)$, $m\left(3\sqrt{2}\right)$, and
$m\left(\ii \sqrt{2}\right)$}\label{sec:sqrt{2}}
Setting
$n=\frac{1}{\sqrt{2}}$ in Eq. \eqref{eq:first}, we obtain
\[m\left(3\sqrt{2}\right)+m\left(\ii \sqrt{2}\right) =m(8).\]
Doing the same in Eq. \eqref{eq:ko}, we find that
\[m(2) + m(8)=2m\left(3\sqrt{2}\right). \]
In this section we will establish the identity
\[3m(3\sqrt{2})=5m(\ii \sqrt{2}),\]
from which we can deduce expressions for $m(2)$ and $m(8)$.

Consider the functions $f$ and $1-f$, where $f=\frac{\sqrt{2} Y
-X}{2} \in \mathbb{C}(E_{3\sqrt{2}})$. Their divisors are
\begin{align*}
\left(\frac{\sqrt{2} Y -X}{2} \right) =&(2P)+2(P+Q)-3O,\\
\left(1-\frac{\sqrt{2} Y -X}{2} \right) =& (P)+(Q)+(3P+Q)-3O.
\end{align*}
The diamond operation yields
\[(f)\diamond(1-f) = 6(P)-10(P+Q).\]
But $(f)\diamond(1-f)$ is trivial in $K$-theory, hence
\[6(P) \sim 10(P+Q).\]

Now consider the isomorphism $\phi$:
\begin{equation} \label{eq:iso}
\phi: E_{2\left(n+\frac{1}{n}\right)} \rightarrow E_{2\left(\ii
n+\frac{1}{\ii n}\right)}, \qquad (X,Y) \rightarrow (-X, \ii Y)
\end{equation}
This isomorphism implies that
\[\r_{\ii \sqrt{2}}\left(\left\{x,y\right\}\right) = \r_{3
\sqrt{2}}\left(\left\{x \circ\phi,y\circ\phi\right\}\right). \]
But we know that
\[\left(x\circ \phi\right) \diamond \left(y\circ \phi\right)  = 8 (P+Q).\]

This implies
\[ 6 \r_{3 \sqrt{2}} \left( \left\{ x,y \right\}\right) = 10 \r_{\ii\sqrt{2}}(\left\{ x,  y\right\}),\]
and
\[3m(3\sqrt{2})=5m(\ii \sqrt{2}).\]
Consequently, we may conclude that
\[m(8) = \frac{8}{5}m(3\sqrt{2}), \quad m(2) = \frac{2}{5}m(3\sqrt{2}),\]
and finally
\[m(8)=4m(2).\]

\subsection{The Hesse family}
We will now sketch the case of the Hesse family:
\[x^3+y^3+1-\frac{3}{t^{\frac{1}{3}}} xy.\]

This family corresponds to $\Gamma_0(3)$. The diamond operation yields
\begin{equation}
(x) \diamond (y)  =9(P)+9(A)+9(B),
\end{equation}
where $P$ is a point of order 3, defined over $\mathbb{Q}(t^{1/3})$, and $A, B$ are points of order 3 such that $A+B+P=O$.

For $0<t<1$, we have
 \[n(t) = \frac{9}{2\pi}  \Im
\left(\frac{\tau}{y_\tau}\left(
R_\tau\left(\e^{\frac{4\pi\ii}{3}}\right)
+R_\tau\left(\e^{\frac{4\pi\ii(1+\tau)}{3}}\right) +
R_\tau\left(\e^{\frac{ 2\pi\ii (2+\tau)}{3}}\right) \right)
\right).\] If we let $\mu = -\frac{1}{\tau}$, we obtain, after
several steps,
\[n(t)=\Re \left(\frac{27 \sqrt{3} y_\mu}{4\pi^2} \sum_{k, n}'
\frac{\chi_{-3}(n) } {(3\mu k+n)^2(3\bar{\mu}k+n)} \right).\]

Following the previous example, this result may be extended to
$\mathbb{C}\setminus \kappa$ by comparing holomorphic functions.


\subsection{The $\Gamma_0^0(6)$ example}
We will now sketch a treatment of Stienstra's example \cite{St}:
\[(x+1)(y+1)(x+y)-\frac{1}{t}xy.\]

Applying the diamond operation, we have
\[(x)\diamond(y)=-6(P)-6(2P),\]
where $P$ is a point of order 6.

For $t$ small, one can write
\[g(t)= \frac{3}{\pi} \Im\left(\frac{\tau}{y_\tau} R_\tau(\xi_6^{-1}) +R_\tau(\xi_3^{-1})\right).\]
Eventually, one arrives to
\begin{equation*}
g(t) = \Re\left( \frac{36y_\mu}{\pi^2} \sum_{m  ,n}'
\frac{\chi_{-3}(m)}{(m+6\mu n)^2(m+6\bar{\mu}n)}\right)
 +
\Re\left( \frac{9y_\mu}{2\pi^2} \sum_{m  ,n}'
\frac{\chi_{-3}(m)}{(m+3\mu n)^2(m+3\bar{\mu}n)}\right),
\end{equation*}
thus recovering a result of \cite{St}.

\subsection{The $\Gamma_0^0(5)$ example}
Now we will consider our final example:
\[(x+y+1)(x+1)(y+1)-\frac{1}{t} xy.\]

Applying the diamond operation, we find that
\[(x)\diamond (y) = 10(P)+5(2P),\]
where $P$ is a torsion point of order 5.

For $t>0$
\[r(t)= \frac{5}{2\pi} \Im \left( \frac{\tau}{y_\tau}
\left(2R_\tau\left(\e^{\frac{8\pi\ii}{5}}\right)
+R_\tau\left(\e^{\frac{6\pi\ii}{5}}\right)\right)\right).\]
Finally,
\[r(t)=-\Re\left( \frac{25\ii y_\mu}{4\pi^2} \sum_{m,n}'
\frac{2\left(\zeta_5^m-\zeta_5^{-m} \right)+
\zeta_5^{2m}-\zeta_5^{-2m} }{(m+5\mu n)^2(m+5\bar{\mu}n)}\right).\]

In conclusion, we see that the modular structure comes from the form of the regulator function, and the functional identities are consequences of the functional identities of the elliptic dilogarithm.

\section{Conclusion}
We have used both regulator and $q$-series methods to prove a
variety of identities between the Mahler measures of genus-one
polynomials. We will conclude this paper with a final open problem.\\
\textbf{Open Problem 3: } How do you characterize all the functional
equations of $\mu(t)$?

We have seen that there are identities like Eq. \eqref{eq:ko},
stating that
\begin{equation*}
\begin{split}
2\m\left(2\left(k+\frac{1}{k}\right)+x+\frac{1}{x}+y+\frac{1}{y}\right)=&\m\left(4k^2+x+\frac{1}{x}+y+\frac{1}{y}\right)\\
&+\m\left(\frac{4}{k^2}+x+\frac{1}{x}+y+\frac{1}{y}\right).
\end{split}
\end{equation*}
While this formula does not follow from Eq. \eqref{Hecke eigenvalue
equation}, it can be proved with regulators.

Indeed, the last section showed us that we can obtain functional
identities for the Mahler measures by looking at functional
equations for the elliptic dilogarithm.

Now, understanding these identities is a very hard problem. To have
an idea of the dimensions of this problem, let us note that equation
\eqref{eq:JNk} corresponds to the integration of an identity for the
Hecke operator $T_p$. This suggests that more identities will follow
from looking at the general operator $T_n$.  And this is just the
beginning of the story\dots


\bigskip
\begin{acknowledgements}

The authors would like to deeply thank  David Boyd and Fernando
Rodriguez-Villegas for many helpful discussions, and David Boyd for pointing out the work of Kurokawa and Ochiai \cite{KO}.  ML extends her gratitude to Christopher
Deninger, Herbert Gangl, and Florian Herzig for enlightening discussions.

ML is a postdoctoral fellow at the Pacific Institute for the Mathematical Sciences and the University of British Columbia. This research was also partially conducted while ML was a member at the Institute for
Advanced Study, and at the Mathematical Sciences Research Institute, a visitor
at the Institut des Hautes \'Etudes Scientifiques,  a guest at the
Max-Planck-Insitut f\"ur Mathematik, and she was employed by the Clay Mathematics Institute as a Liftoff Fellow. ML whishes to thank these institutions for
their support and hospitality.

This material is partially based upon work supported by the National Science Foundation under agreement No. DMS-0111298.

\end{acknowledgements}


\end{document}